\documentclass[12pt, oneside]{article}
\usepackage{amsthm}
\usepackage{amsthm}
\usepackage{graphicx} 
\usepackage{tikz}
\usepackage{mathrsfs}
\usepackage{fancyhdr}
\usetikzlibrary{decorations.fractals}
\usepackage[english]{babel}
\usepackage[many]{tcolorbox}
\usepackage{geometry}
\usepackage{amsfonts,amsmath}
\usepackage{amsmath}
\usepackage{tkz-euclide}
\usepackage{hyperref}
\hypersetup{colorlinks=true,linkcolor=blue,filecolor=magenta,urlcolor=cyan}
\usepackage{float}
\usepackage{lscape}
\usepackage{setspace}
\newtheorem{ex}{Example}[subsection]
\newtheorem{Theorem}{Theorem}[subsection]
\newtheorem{Lemma}{Lemma}[subsection]

\usepackage{mathtools}

\title{\textbf{Some Classes of series involving the Riemann zeta function, Fibonacci numbers and the Lucas numbers}}
\author{Akerele Olofin Segun
	\\ Department of Mathematics, University of Ibadan, Oyo State, Nigeria.\\
	email: akereleolofin@gmail.com}
\date{}
\AtBeginDocument{}

\begin{document}
	\maketitle
		\fancyhead{}
	\begin{abstract}
		The objective of this manuscript is to offer explicit expressions for diverse categories of infinite series incorporating the Fibonacci (Lucas) sequence and the Riemann zeta function. In demonstrating our findings, we will utilize conventional methodologies and integrate the Binet formulas pertinent to these sequences with generating functions that encompass the Riemann zeta function alongside established evaluations of certain series.
	\end{abstract}
	\section{Introduction}
	\par This manuscript is dedicated to the amalgamation of two profoundly eminent and pivotal mathematical entities: the Riemann zeta function and Fibonacci numbers. Despite the extensive scrutiny and comprehension of both constructs individually, it is notably intriguing that documented correlations between them are conspicuously absent within the annals of mathematical literature. This anomaly serves as an impetus for our scholarly endeavor. Within this exposition, our objective is to rectify this lacuna. Meanwhile, this paper was motivated by an article, General Infinite series evaluations involving Fibonacci numbers and the Riemann zeta function by Robert.F and Taras.G.\cite{3}
	. \par Specifically, we embark on an exhaustive examination of various categories of infinite series that incorporate these illustrious mathematical constructs. Through the meticulous utilization of generating functions in conjunction with methodical reasoning, we shall endeavor to ascertain closed-form representations for these series. It is noteworthy that the majority of these infinite series, in their manifestation, exhibit simplistic analytical expressions, frequently characterized by the involvement of elementary trigonometric functions.
	\par Recall that the polygamma function $\psi^{(m)}(z)$ of order $m$ is a meromorphic function on the complex numbers $\mathbb{C}$ defined as the ($m+1$)th derivative of the logarithm of the gamma function $\Gamma(z)$ \cite{1}:
	\begin{equation*}
		\psi^{(m)}(z):= \frac{d^m}{dz^m}\psi(z)=\frac{d^{m+1}}{dz^{m+1}}\ln \Gamma(z).
	\end{equation*}
	Thus, $\psi^{(0)}(z)=\psi(z)=\displaystyle{\frac{\Gamma '(z)}{\Gamma(z)}}$ holds where  $\psi^{(0)}(z)$ is the digamma function, meanwhile they are holomorphic on $\mathbb{C}\setminus \mathbb{Z}_{\leq 0}$.\\
	The polygamma function possesses the following properties: 
	\begin{equation}
		\psi^{(m)}(z+1)=\psi^{(m)}(z)+\frac{(-1)^m m!}{z^{m+1}},
	\end{equation}
	\vspace{-5mm}
	\begin{equation}
		\label{2}
		\psi^{(m)}(z)=(-1)^{m+1}m!\sum_{k=0}^{\infty}\frac{1}{(z+k)^{m+1}}, \quad m>0,z\in \mathbb{C}\setminus \mathbb{Z}_{\leq 0}
	\end{equation}
	
	And,  $\psi^{(0)}(n)=-\gamma + \displaystyle{\sum_{k=1}^{n-1}\frac{1}{k}}$ for all $n\in \mathbb{N}$. \\
	We also have the reflection property given as; $\displaystyle{(-1)^m \psi^{(m)}(1-z)- \psi^{(m)}(z)=\pi \frac{d^m}{dz^m}\cot(\pi z)}$, where $\gamma$ is the Euler-Mascheroni constant given by
	$$\gamma = \lim_{n\to \infty}\left(\sum_{k=1}^n\frac{1}{k}-\ln n \right)=0.5772156649\dots $$
	\par On the other hand, the Riemann zeta function $\zeta(s)$, $s\in \mathbb{C}$, is defined by \cite{1}
	$$\zeta(s)=\sum_{k=1}^{\infty}\frac{1}{k^s}, \quad   \Re(s)>1.$$
	The analytical continuation to all $s\in \mathbb{C}$ with $\Re(s)>0, s\neq 1$, is given by $$\zeta(s)=\left(1-2^{1-s}\right)^{-1}\sum_{k=1}^{\infty}\frac{(-1)^{k+1}}{k^s}.$$
	\vspace{-7mm}
	\par Meanwhile, let $F_n$ and $L_n$ denote the n-th Fibonacci and Lucas numbers, both satisfying the recurrence relation $\Gamma_n=\Gamma_{n-1}+\Gamma_{n-2}$, $n\geq 2$, with conditions $F_0=0,F_1=1$ and $L_0=2,L_1=1$. Also $L_{-m}=(-1)^mL_m$ and $F_{-m}=(-1)^{m-1}F_m$. Through out this paper, we denote the golden ratio $\alpha=\frac{1+\sqrt{5}}{2}$ and write it's conjugate $\beta=\frac{1-\sqrt{5}}{2}$, so that $\alpha \beta = -1$ and $\alpha + \beta = 1$. We have the Binet formulas for Fibonacci and Lucas numbers to be; 
	$$F_m = \frac{\alpha^m-\beta^m}{\alpha-\beta}, \quad L_m = \alpha^m + \beta^m$$
	for any integer $m$.\\
	The Fibonacci and Lucas sequence are indexed in the On-Line Encyclopedia of Integer Sequences \cite{7} with entries A000045 and A000032, respectively.
	\par It is known that any real number has representation of the form $\sum_{k=2}^{\infty}q_k(\zeta(k)-1)$, where the rational coefficients $q_k$ are, in some appropriate sense, well behaved, (See Pg 261, \cite{2}). We focus on the cases where $q_k$ is a function of $F_k$ or $L_k$. Among the conclusions drawn in this paper, we will ascertain that for $n\geq 2$ and $n\geq 1$ respectively (see \S \ \ref{some} ) $$\sum_{k=1}^{\infty}\frac{\zeta(2k)}{kF^{4k}_{2n}}=\ln\left(\frac{\pi}{F_{2n}^2}\right)+\ln\csc\frac{\pi}{F^2_{2n}} , \ \sum_{k=2}^{\infty}\frac{\zeta(k)}{k(F_nL_n)^k}=\ln \Gamma\left(\frac{F_{2n}-1}{F_{2n}}\right)-\frac{\gamma}{F_{2n}}$$
	\par Evaluations of infinite series that incorporate both Fibonacci numbers and the zeta function are seldom encountered.Thus, we have some interesting identities involving $\zeta(s)$ and Fibonacci (Lucas) numbers stated in (Pg 119-120, \cite{4}): 
	$$\sum_{k=1}^{\infty}\zeta(2k+1)\frac{F_{2k}}{5^k}=\sum_{k=1}^{\infty}(\zeta(2k+1)-1)F_{2k}=\frac{1}{2}, \quad \sum_{k=2}^{\infty}(\zeta(k)-1)F_{k-1}=1+\frac{\pi}{\sqrt{5}}\tan\frac{\sqrt{5}\pi}{2}.$$
	\par The objective of this article is to further research in this area and to provide additional closed-form expressions for certain categories of infinite series that incorporate Fibonacci (Lucas) numbers and the Riemann zeta function. Throughout this paper, we verify our results using the Computer Algebra System (CAS) software \textsf{Mathematica 13.3}.
	\section{Main Results}
	\begin{Theorem}
		For $m\geq 2$, we have
		\begin{equation}
			\label{3}
			\sum_{k=1}^{\infty}(2k-1)(\zeta(2k)-1)F_{2k+m+1}=\frac{\pi^2}{2}\sec^2\left(\frac{\sqrt{5}\pi}{2}\right)F_{m+3}-\frac{9F_m+7F_{m-1}}{2},
		\end{equation}
		\vspace{-5mm}
		\begin{equation}
			\label{4}
			\sum_{k=1}^{\infty}(2k-1)(\zeta(2k)-1)L_{2k+m+1}=\frac{\pi^2}{2}\sec^2\left(\frac{\sqrt{5}\pi}{2}\right)L_{m+3}-\frac{9L_m+7L_{m-1}}{2}.
		\end{equation}
	\end{Theorem}
	\begin{proof}
		From (See Pg 281, \cite{6}) we have that 
		$$\sum_{k=1}^{\infty}(\zeta(2k)-1)z^{2k-1}=-\frac{\pi}{2}\cot \pi z + \frac{3z^2-1}{2z(z^2-1)}, \quad |z|<1$$ On differentiating and multiplying by $z^2$, we have 
		$$\sum_{k=1}^{\infty}(2k-1)(\zeta(2k)-1)z^{2k}=\frac{1}{2}\left(\pi^2 z^2\csc^2 \pi z -\frac{3z^4+1}{(z^2-1)^2}\right).$$ Next we set $z=\alpha$ while using $\alpha^2=\alpha+1$ and $\alpha^4=3\alpha+2$ we get that
		$$\sum_{k=1}^{\infty}(2k-1)(\zeta(2k)-1)\alpha^{2k}=\frac{1}{2}\left(\pi^2 \alpha^2\csc^2 \pi \alpha -\frac{9\alpha+7}{\alpha^2}\right).$$
		It is clear that for $m\geq 2$, it follows that 
		\begin{equation}
			\label{5}
			\sum_{k=1}^{\infty}(2k-1)(\zeta(2k)-1)\alpha^{2k+m+1}=\frac{1}{2}\left(\pi^2 \alpha^{m+3}\csc^2 \pi \alpha -\frac{9\alpha+7}{\alpha^{1-m}}\right).
		\end{equation}
		Similarly, we get 
		\begin{equation}
			\label{6}
			\sum_{k=1}^{\infty}(2k-1)(\zeta(2k)-1)\beta^{2k+m+1}=\frac{1}{2}\left(\pi^2 \beta^{m+3}\csc^2 \pi \beta -\frac{9\beta+7}{\beta^{1-m}}\right).
		\end{equation}
		Now joining (\ref{5}) and (\ref{6}) using the Binet formula we have,
		$$\sum_{k=1}^{\infty}(2k-1)(\zeta(2k)-1)L_{2k+m+1}=\frac{\pi^2}{2}\left(\alpha^{m+3}\csc^2 \pi \alpha + \beta^{m+3}\csc^2 \pi \beta \right)-\frac{1}{2}\left(\frac{9\alpha+7}{\alpha^{1-m}}+\frac{9\beta+7}{\beta^{1-m}}\right)$$
		$$= \frac{\pi^2}{2}\left(\alpha^{m+3}\sec^2\left(\frac{\sqrt{5}\pi}{2}\right) + \beta^{m+3}\sec^2\left(\frac{\sqrt{5}\pi}{2}\right) \right)-\frac{1}{2}(9L_m+7L_{m-1})$$
		$$=\frac{\pi^2}{2}\sec^2\left(\frac{\sqrt{5}\pi}{2}\right)L_{m+3}-\frac{9L_m+7L_{m-1}}{2}.$$
		As desired. Similarly we can prove (\ref{3}) but the proof is omitted and left to the reader.
	\end{proof}
	\begin{ex}
		For $m=2$ and $m=3$ we have the following,
		\begin{align*}
			\sum_{k=1}^{\infty}(2k-1)(\zeta(2k)-1)F_{2k+3} &=\frac{5\pi^2}{2}\sec^2\left(\frac{\sqrt{5}\pi}{2}\right)-8,\\
			\sum_{k=1}^{\infty}(2k-1)(\zeta(2k)-1)F_{2k+4} &= 4\pi^2\sec^2\left(\frac{\sqrt{5}\pi}{2}\right)-\frac{25}{2},\\
			\sum_{k=1}^{\infty}(2k-1)(\zeta(2k)-1)L_{2k+3} &=\frac{11\pi^2}{2}\sec^2\left(\frac{\sqrt{5}\pi}{2}\right)-17,\\
			\sum_{k=1}^{\infty}(2k-1)(\zeta(2k)-1)L_{2k+4} &=9\pi^2\sec^2\left(\frac{\sqrt{5}\pi}{2}\right)-\frac{57}{2}.
		\end{align*}
	\end{ex}
	\begin{Theorem}
		For $m \geq 0$, we have
		\begin{multline*}
			\sum_{k=2}^{\infty}\frac{(k+1)(\zeta(k)-1)}{k}L_{k+m} 
			\vspace{-5mm} =2(1-\gamma)L_{m+1}+\alpha^m\left(\ln \Gamma\left(\beta^2\right)-\alpha\psi^{(0)}\left(\beta^2\right)\right)\\ + 
			\beta^m\left(\ln \Gamma\left(\alpha^2\right) -\beta\psi^{(0)}\left(\alpha^2\right)\right)
		\end{multline*}
		\vspace{-8mm}
		\begin{multline*}
			\sum_{k=2}^{\infty}\frac{(k+1)(\zeta(k)-1)}{k}F_{k+m}  
			\vspace{-3mm} =2(1-\gamma)F_{m+1}+\frac{\alpha^m}{\sqrt{5}}\left(\ln \Gamma\left(\beta^2\right)-\alpha\psi^{(0)}\left(\beta^2\right)\right) \\ - \frac{\beta^m}{\sqrt{5}}\left(\ln \Gamma\left(\alpha^2\right)-\beta\psi^{(0)}\left(\alpha^2\right)\right)
		\end{multline*}
	\end{Theorem}
	\begin{proof}
		We can show using (See Pg 280, \cite{6}), that 
		$$\sum_{k=2}^{\infty}(\zeta(k)-1)\frac{z^k}{k}=(1-\gamma)z+\ln\Gamma(2-z), \quad |z|<2$$
		By multiplying the above by $z$ and then differentiating, it follows that,
		\begin{equation}
			\label{7}
			\sum_{k=2}^{\infty}\frac{(k+1)(\zeta(k)-1)}{k}z^k=2(1-\gamma)z+\ln \Gamma(2-z)-z\psi^{(0)}(2-z), \quad |z|<2
		\end{equation}
		Now, set $z=\alpha$ then we have 
		$$	\sum_{k=2}^{\infty}\frac{(k+1)(\zeta(k)-1)}{k}\alpha^k=2(1-\gamma)\alpha+\ln \Gamma\left(\beta^2\right)-\alpha\psi^{(0)}\left(\beta^2\right)$$
		Thus, 
		$$	\sum_{k=2}^{\infty}\frac{(k+1)(\zeta(k)-1)}{k}\left(\alpha^{k+m}+\beta^{k+m}\right)$$
		$$=2(1-\gamma)\left(\alpha^{m+1}+\beta^{m+1}\right)+\beta^m\ln \Gamma\left(\alpha^2\right)+\alpha^m \ln \Gamma\left(\beta^2\right)-\left(\alpha^{m+1}\psi^{(0)}\left(\beta^2\right)+\beta^{m+1}\psi^{(0)}\left(\alpha^2\right)\right)$$
		By using the Binet formula we get that,
		$$
		\sum_{k=2}^{\infty}\frac{(k+1)(\zeta(k)-1)}{k}L_{k+m} $$
		\vspace{-5mm} $$=2(1-\gamma)L_{m+1}+\alpha^m\left(\ln \Gamma\left(\beta^2\right)-\alpha\psi^{(0)}\left(\beta^2\right)\right)\\+ 
		\beta^m\left(\ln \Gamma\left(\alpha^2\right)-\beta\psi^{(0)}\left(\alpha^2\right)\right)$$
		The preceding identity can be proved in similarly.
	\end{proof}
	
	At this point we want to investigate $\alpha^{m+1}\psi^{(0)}\left(\beta^2\right)+\beta^{m+1}\psi^{(0)}\left(\alpha^2\right)$ and \\ $\frac{1}{\sqrt{5}}\left(\beta^{m+1}\psi^{(0)}\left(\alpha^2\right)-\alpha^{m+1}\psi^{(0)}\left(\beta^2\right)\right)$.
	Before we proceed we recall 
	from (\ref{2}),
	$$\psi^{(0)}\left(\beta^2\right)=-\gamma + \sum_{n=1}^{\infty}\left(\frac{1}{n}-\frac{1}{n+\beta}\right) \ and \  \psi^{(0)}\left(\alpha^2\right)=-\gamma + \sum_{n=1}^{\infty}\left(\frac{1}{n}-\frac{1}{n+\alpha}\right). $$
	It follows that,
	$$\alpha^{m+1}\psi^{(0)}\left(\beta^2\right)+\beta^{m+1}\psi^{(0)}\left(\alpha^2\right)=-\gamma \left(\alpha^{m+1}+\beta^{m+1}\right)+\sum_{n=1}^{\infty}\left(\frac{\alpha^{m+1}}{n}-\frac{\alpha^{m+1}}{n+\beta}+\frac{\beta^{m+1}}{n}-\frac{\beta^{m+1}}{n+\alpha}\right)$$
	$$=-\gamma L_{m+1}+\sum_{n=1}^{\infty}\left(\frac{(n-1)L_{m+1}-nL_{m+2}}{n(n+\alpha)(n+\beta)}\right)=-\gamma L_{m+1}+\sum_{n=1}^{\infty}\left(\frac{n L_{m+1}-L_{m+1}-nL_{m+2}}{n(n^2+n-1)}\right)$$
	$$-\gamma L_{m+1}+\sum_{n=1}^{\infty}\frac{L_{m+1}}{n^2+n-1}-\sum_{n=1}^{\infty}\frac{L_{m+1}}{n(n^2+n-1)}-\sum_{n=1}^{\infty}\frac{L_{m+2}}{n^2+n-1}$$
	$$=-\gamma L_{m+1}+\underbrace{\left(L_{m+1}-L_{m+2}\right)}_{-L_m}\sum_{n=1}^{\infty}\frac{1}{n^2+n-1}-L_{m+1}\sum_{n=1}^{\infty}\frac{1}{n(n^2+n-1)}$$
	The series above was already evaluated in \cite{5}, Note that 
	$$\sum_{n=1}^{\infty}\frac{1}{n^2+n-1}=1+\frac{\sqrt{5}\pi}{5}\tan\frac{\sqrt{5}\pi}{2} , \ \sum_{n=1}^{\infty}\frac{1}{n(n^2+n-1)}=1-\gamma -\psi^{(0)}(\alpha)+\frac{5+\sqrt{5}}{10}\pi \tan\frac{\sqrt{5}\pi}{2} .$$
	Thus, on simplification we obtain, 
	$$
	\alpha^{m+1}\psi^{(0)}\left(\beta^2\right)+\beta^{m+1}\psi^{(0)}\left(\alpha^2\right)$$ $$=\left(\psi^{(0)}(\alpha)-1-\frac{5+\sqrt{5}}{10}\pi \tan\frac{\sqrt{5}\pi}{2}\right)L_{m+1}-\left(1+\frac{\sqrt{5}\pi}{5}\tan\frac{\sqrt{5}\pi}{2}\right)L_m$$
	Similarly, 
	\begin{multline*}
		\frac{1}{\sqrt{5}}\left(\beta^{m+1}\psi^{(0)}\left(\alpha^2\right)-\alpha^{m+1}\psi^{(0)}\left(\beta^2\right)\right)\\ =\left(\frac{-\beta^{m+1}}{\sqrt{5}}+\frac{\alpha^{m+1}}{\sqrt{5}}\right)\gamma  +  \frac{1}{\sqrt{5}}\sum_{n=1}^{\infty}\left(\frac{\beta^{m+1}}{n}-\frac{\beta^{m+1}}{n+\alpha}-\frac{\alpha^{m+1}}{n}+\frac{\alpha^{m+1}}{n+\beta}\right)\\ = \gamma F_{m+1}+\frac{1}{\sqrt{5}}\sum_{n=1}^{\infty}\left(\frac{(1-n)\sqrt{5}F_{m+1}+n\sqrt{5}F_{m+2}}{n(n+\alpha)(n+\beta)}\right)
		\\ =\gamma F_{m+1}+\sum_{n=1}^{\infty}\frac{F_{m+1}}{n(n^2+n-1)}-\sum_{n=1}^{\infty}\frac{F_{m+1}}{n^2+n-1}+\sum_{n=1}^{\infty}\frac{F_{m+2}}{n^2+n-1}
		\\=\gamma F_{m+1}+\underbrace{\left(F_{m+2}-F_{m+1}\right)}_{F_m}\sum_{n=1}^{\infty}\frac{1}{n^2+n-1}+F_{m+1}\sum_{n=1}^{\infty}\frac{1}{n(n^2+n-1)}
		\\=F_{m+2}+\left(\frac{\sqrt{5}\pi}{5}\tan\frac{\sqrt{5}\pi}{2}\right)F_m+\left(\frac{5+\sqrt{5}}{10}\pi \tan\frac{\sqrt{5}\pi}{2}-\psi^{(0)}(\alpha)\right)F_{m+1}
	\end{multline*}
	With these, we have another strong closed form for Theorem 2.0.2.
	For m=0 we have the following,
	$$\sum_{k=2}^{\infty}\frac{(k+1)(\zeta(k)-1)}{k}L_k = 5-2\gamma-\psi^{(0)}\left(\frac{1+\sqrt{5}}{2}\right)+\ln\left(-\pi \sec \frac{\sqrt{5}\pi}{2}\right)+\frac{1+\sqrt{5}}{2}\pi\tan\frac{\sqrt{5}\pi}{2}$$
	$$	\sum_{k=2}^{\infty}\frac{(k+1)(\zeta(k)-1)}{k}F_{k}=3-2\gamma +\frac{1}{\sqrt{5}}\ln \left[\frac{\Gamma\left(\frac{3-\sqrt{5}}{2}\right)}{\Gamma\left(\frac{3+\sqrt{5}}{2}\right)}\right]+\frac{5+\sqrt{5}}{10}\pi\tan\frac{\sqrt{5}\pi}{2}-\psi^{(0)}\left(\frac{1+\sqrt{5}}{2}\right) $$ 
	\vspace{2mm}
	\begin{Theorem}
		If $\omega (m,x)=x^m\left(\psi^{(0)}(x+2)+x\psi^{(1)}(x+2)+\psi^{(0)}(2-x)-x\psi^{(1)}(2-x)\right)$ for $x\in (-2,2)$, we have 
		\begin{equation}
			\sum_{k=1}^{\infty}(2k+1)(\zeta(2k+1)-1)F_{2k+m} =  (1-\gamma)F_m - \frac{1}{2\sqrt{5}}(\omega(m,\alpha)-\omega(m,\beta)), \quad m\geq 0
		\end{equation}
		\vspace{-5mm}
		\begin{equation}
			\sum_{k=1}^{\infty}(2k+1)(\zeta(2k+1)-1)L_{2k+m} = (1-\gamma)L_m - \frac{1}{2}(\omega(m,\alpha)+\omega(m,\beta)), \quad m\geq 0
		\end{equation}
	\end{Theorem}
	\begin{proof}
		We can show that, 
		\begin{equation}
			\label{10}
			\sum_{k=1}^{\infty}(\zeta(2k+1)-1)z^{2k} = (1-\gamma)-\frac{1}{2}\left(\psi^{(0)}(2+z)+\psi^{(0)}(2-z)\right), \quad |z|<2
		\end{equation}
		But the proof is left to the reader or check (Pg 280, \cite{6}), Now we multiply (\ref{10}) by z and differentiate to get,
		\begin{multline*}
			\sum_{k=1}^{\infty}(2k+1)(\zeta(2k+1)-1)z^{2k} = (1-\gamma)-\frac{d}{dz}\left(\frac{z}{2}\left(\psi^{(0)}(2+z)+\psi^{(0)}(2-z)\right)\right)\\ = (1-\gamma)-\frac{1}{2}\left(\psi^{(0)}(z+2)+z\psi^{(1)}(z+2)+\psi^{(0)}(2-z)-z\psi^{(1)}(2-z)\right), 
		\end{multline*}
		The remainder of the proof is as above and we leave it as an exercise.
	\end{proof}
	\begin{ex}
		When $m=0$, we see that the following evaluations are involving the Riemann zeta function at positive odd integer argument and scaled even Fibonacci (Lucas) numbers,
		\begin{align*}
			\sum_{k=1}^{\infty}(2k+1)(\zeta(2k+1)-1)F_{2k}&=\frac{1}{2\sqrt{5}}(\omega(0,\beta)-\omega(0,\alpha)),\\
			\sum_{k=1}^{\infty}(2k+1)(\zeta(2k+1)-1)L_{2k}&=2(1-\gamma)-\frac{1}{2}\left(\omega(0,\alpha)+\omega(0,\beta)\right).
		\end{align*}
	\end{ex}
	\section{Some Class of series involving the Riemann zeta function and reciprocals of the Fibonacci (Lucas) Numbers}\label{some}
	Before we proceed, we recall a more general form of the Riemann zeta function, popularly known as the Hurwitz zeta function $\zeta(s,a)$, $s\in \mathbb{C}$ defined by \cite{1}
	$$\zeta(s,a)=\sum_{n=0}^{\infty}\frac{1}{(n+a)^s}, \quad a \neq 0,-1,-2,\dots $$
	where $\Re(s)>1$. This series is absolutely convergent for the given values of $s$ and $a$ and can be extended to a meromorphic function defined for all $s\neq 1$. The Riemann zeta function is of course $\zeta(s,1)$.
	\par Now, we present several lemmas essential for this series category.
	\subsection{Lemmas}
	\begin{Lemma}
		\begin{equation}
			\label{11}
			\sum_{k=2}^{\infty}\zeta(k,a)\frac{z^k}{k}=\ln \Gamma(a-z)-\ln \Gamma(a)+z\psi^{(0)}(a), \quad (|z|<|a|)
		\end{equation}
		\begin{proof}
			It can be shown using (Pg 276, \cite{6}) that,
			$$\sum_{k=2}^{\infty}(-1)^k\zeta(k,a)\frac{t^k}{k}=\ln \Gamma(a+t)-\ln \Gamma(a)-t\psi^{(0)}(a), \quad (|t|<|a|)$$
			Setting $t=-z$ concludes the proof.
		\end{proof}
	\end{Lemma}
	\begin{Lemma}
		\begin{equation}
			\label{12}
			\sum_{k=1}^{\infty}\zeta(2k,a)\frac{z^{2k}}{k}=\ln \Gamma(a+z)+\ln \Gamma(a-z)-2\ln \Gamma(a), \quad (|z|<|a|)
		\end{equation}
	\end{Lemma}
	\begin{proof}
		This is just a direct consequence of (\ref{11}), the proof is left to interested readers.
	\end{proof}
	\begin{Lemma}
		\begin{equation}
			\label{13}
			\sum_{k=1}^{\infty}\zeta(2k+1,a)\frac{z^{2k+1}}{2k+1}=\frac{1}{2}\left(\ln \Gamma(a-z)-\ln \Gamma (a+z)\right)+z\psi^{(0)}(a), \quad (|z|<|a|)
		\end{equation}
	\end{Lemma}
	\begin{proof}
		By differentiating both sides of (\ref{12}) with respect to $z$ and then multiply by $z$, the result is straightforward to deduce.
	\end{proof}
	\subsection{Results}
	\begin{Theorem}
		Let $m$ and $n$ be positive integers such that $n\geq m$, unless stated otherwise. Then
		\begin{align*}
			\sum_{k=2}^{\infty}\frac{\zeta(k)}{k(L_nF_m+F_nL_m)^k}&=\ln\Gamma\left(\frac{2F_{n+m}-1}{2F_{n+m}}\right)-\frac{\gamma}{2F_{n+m}},\\
			\sum_{k=2}^{\infty}\frac{\zeta(k)}{k(F_n^2+(-1)^{n+m-1}F_m^2)^k}&=\ln\Gamma\left(\frac{F_{n-m}F_{n+m}-1}{F_{n-m}F_{n+m}}\right)-\frac{\gamma}{F_{n-m}F_{n+m}},\quad n>m\\
			\sum_{k=2}^{\infty}\frac{\zeta(k)}{k(L_{n+m}+(-1)^mL_{n-m})^k}&=\ln\Gamma\left(\frac{L_mL_n-1}{L_mL_n}\right)-\frac{\gamma}{L_mL_n},
		\end{align*}
		\vspace{-3mm}
		$$\sum_{k=2}^{\infty}\frac{\zeta(k)}{k(5F_mF_n)^k}= \ln \Gamma\left(\frac{L_{n+m}+(-1)^{m-1}L_{n-m}-1}{L_{n+m}+(-1)^{m-1}L_{n-m}}\right)-\frac{\gamma}{L_{n+m}+(-1)^{m-1}L_{n-m}},$$
		$$\sum_{k=2}^{\infty}\frac{\zeta(k)}{k(F_mL_n)^k}= \ln \Gamma\left(\frac{F_{n+m}+(-1)^{m-1}F_{n-m}-1}{F_{n+m}+(-1)^{m-1}F_{n-m}}\right)-\frac{\gamma}{F_{n+m}+(-1)^{m-1}F_{n-m}}.$$
	\end{Theorem}
	\begin{proof}
		To show these results, we require the following well-known identities \cite{4}:
		\begin{equation}
			\label{14}
			F^2_n+(-1)^{n+m-1}F^2_m=F_{n-m}F_{n+m},
		\end{equation}
		\vspace{-10mm}
		\begin{equation}
			\label{15}
			F_{n+m}+(-1)^mF_{n-m}=L_mF_n,
		\end{equation}
		\vspace{-8mm}
		\begin{equation}
			\label{16}
			L_{n}F_m+F_nL_m=2F_{n+m},
		\end{equation}
		\vspace{-8mm}
		\begin{equation}
			\label{17}
			L_{n+m}+(-1)^mL_{n-m}=L_mL_n,
		\end{equation}
		\vspace{-8mm}
		\begin{equation}
			\label{18}
			L_{n+m}+(-1)^{m-1}L_{n-m}=5F_mF_n.
		\end{equation}
		Now setting $a=1$ in (\ref{11}) and using (\ref{14})-(\ref{18}), the results follow directly.
	\end{proof}
	\begin{ex}
		The last two identities in Theorem 3.2.1 yields;
		\begin{equation}
			\sum_{k=2}^{\infty}\frac{\zeta(k)}{k(5F^2_n)^k}=\ln\Gamma\left(\frac{L_{2n}+2(-1)^{n-1}-1}{L_{2n}+2(-1)^{n-1}}\right)-\frac{\gamma}{L_{2n}+2(-1)^{n-1}}, \quad n\geq 2
		\end{equation}
		\vspace{-5mm}
		\begin{equation}
			\sum_{k=2}^{\infty}\frac{\zeta(k)}{k(F_nL_n)^k}=\ln \Gamma\left(\frac{F_{2n}-1}{F_{2n}}\right)-\frac{\gamma}{F_{2n}}, \quad n\geq 2
		\end{equation}
	\end{ex}
	\begin{Theorem}
		Let $m$ and $n$ be positive integers such that $n\geq m$, unless stated otherwise. Then
		\begin{align*}
			\sum_{k=1}^{\infty}\frac{\zeta(2k)}{k(F_{n-m}^2F^2_{n+m})^k}&=\ln \left(\frac{\pi}{F^2_{n}+(-1)^{n+m-1}F_m^2}\right)+\ln \csc\left(\frac{\pi}{F^2_{n}+(-1)^{n+m-1}F_m^2}\right), \quad n>m\\
			\sum_{k=1}^{\infty}\frac{\zeta(2k)}{k(L_{m}^2F^2_{n})^k}&=\ln \left(\frac{\pi}{F_{n+m}+(-1)^{m}F_{n-m}}\right)+\ln \csc\left(\frac{\pi}{F_{n+m}+(-1)^{m}F_{n-m}}\right),\\
			\sum_{k=1}^{\infty}\frac{\zeta(2k)}{4^k k F^{2k}_{n+m}}&=\ln\left(\frac{\pi}{L_nF_m+F_nL_m}\right)+\ln\csc\left(\frac{\pi}{L_nF_m+F_nL_m}\right),\\
			\sum_{k=1}^{\infty}\frac{\zeta(2k)}{k(L^2_mL^2_n)^k}&=\ln\left(\frac{\pi}{L_{n+m}+(-1)^{m}L_{n-m}}\right)+\ln\csc\left(\frac{\pi}{L_{n+m}+(-1)^{m}L_{n-m}}\right),\\
			\sum_{k=1}^{\infty}\frac{\zeta(2k)}{25^kk(F^2_mF^2_n)^k}&=\ln\left(\frac{\pi}{L_{n+m}+(-1)^{m-1}L_{n-m}}\right)+\ln\csc\left(\frac{\pi}{L_{n+m}+(-1)^{m-1}L_{n-m}}\right).
		\end{align*}
	\end{Theorem}
	\begin{proof}
		Setting $a=1$ and employing this identity $\Gamma(1+z)\Gamma(1-z)=\pi z \csc\pi z$, in (\ref{12}) and using (\ref{14})-(\ref{18}) then the results follow directly.
	\end{proof}
	\begin{ex}
		Here we have some exquisite identities deduced from Theorem 3.2.2. 
		\begin{equation}
			\sum_{k=1}^{\infty}\frac{\zeta(2k)}{kF^{4k}_{2n}}=\ln\left(\frac{\pi}{F_{2n}^2}\right)+\ln\csc\frac{\pi}{F^2_{2n}}, \quad n\geq 2
		\end{equation}
		\vspace{-5mm}
		\begin{equation}
			\sum_{k=1}^{\infty}\frac{\zeta(2k)}{4^kkF^{2k}_{2n}}=\ln\left(\frac{\pi}{2L_nF_n}\right)+\ln\csc\frac{\pi}{2L_nF_n}, \quad n>0
		\end{equation}
		\vspace{-1mm}
		\begin{equation}
			\sum_{k=1}^{\infty}\frac{\zeta(2k)}{k(L_{2n}F_{2n})^{2k}}=\ln\left(\frac{\pi}{F_{4n}}\right)+\ln\csc\frac{\pi}{F_{4n}}, \quad n>0
		\end{equation}
	\end{ex}
	\begin{Theorem}
		For $n\in \mathbb{N}$. Then
		$$\sum_{k=1}^{\infty}\frac{\zeta(2k+1)}{(2k+1)4^kF^{2k}_{2n}}=F_{2n}\left(\ln \Gamma\left(\frac{2L_nF_n-1}{2L_nF_n}\right)-\ln \Gamma\left(\frac{2L_nF_n+1}{2L_nF_n}\right)\right)-\gamma $$
	\end{Theorem}
	\begin{proof}
		The proof is similar to the previous theorem (See Theorem 3.2.1), thus left to the reader.
	\end{proof}
	\vspace{-4mm}
	Notice that from Theorem 3.2.3, we have this beautiful sum
	\begin{equation}
		\sum_{k=1}^{\infty}\frac{\zeta(2k+1)}{4^k(2k+1)}=\ln 2-\gamma
	\end{equation}
	\section{Conclusion}
	In this article, we presented pioneering closed-form formulations for distinct infinite series that integrate the esteemed Fibonacci and Lucas numbers with the esteemed Riemann zeta function, particularly when its arguments are integers. To authenticate the integrity of our discoveries, we utilized a multifaceted approach, blending Binets formulas with generating functions and established assessments of recognized series.\\ Moreover, through the application of analogous methodologies, we are able to broaden the scope of our findings, extending them to encompass series assessments that establish connections between the Riemann zeta function (with integer arguments) and a spectrum of mathematical entities including Fibonacci and Lucas polynomials, as well as various other renowned sequences of numbers and polynomials.


\begin{thebibliography}{99}
		\bibitem{1} Milton Abramowitz and Irene A Stegun. Handbook of mathematical functions with formulas, graphs, and mathematical tables, volume 55. US Government printing office, 1948.
		\bibitem{2} Jonathan M Borwein, David M Bradley, and Richard E Crandall. Computational strategies for the riemann zeta function. Journal of Computational and Applied Mathematics, 121(1-2):247-296, 2000.
		\bibitem{3} Robert Frontczak and Taras Goy. General infinite series evaluations involving fibonacci numbers and the riemann zeta function. arXiv preprint arXiv:2005.02685, 2020.
		\bibitem{4}Thomas Koshy. Fibonacci and Lucas Numbers with Applications, Volume 2. John Wiley and Sons, 2019.
		\bibitem{5} Anatoli Prudnikov. Integrals and series: special functions.
		\bibitem{6} A Schubert. Et whittaker and gn watson, a course of modern analysis. an introduction to the general theory of infinite processes and of analytic functions; with an account of the principal transcendental functions. 608 s. cambridge 1962. cambridge university press. preis brosch. 27/6 net, 1963.
		\bibitem{7} Neil JA Sloane et al. The on-line encylopedia of integer sequences, 2010.
		\bibitem{8} Hari M Srivastava and Junesang Choi. Zeta and q-Zeta functions and associated series and integrals. Elsevier, 2011.
	\end{thebibliography}
\end{document}